\documentclass{amsart}
\usepackage{amsmath} 
\usepackage{amssymb}
\usepackage{mathrsfs}

\newtheorem{theorem}{Theorem}[section] 
\newtheorem{claim}[theorem]{Claim}

\newtheorem{conclusion}[theorem]{Conclusion}
\newtheorem{observation}[theorem]{Observation}

\theoremstyle{definition}
\newtheorem{definition}[theorem]{Definition}

\newtheorem{explanation}[theorem]{Explanation}
\newtheorem{convention}[theorem]{Convention}

\newtheorem{fact}[theorem]{Fact}

\newtheorem{hypothesis}[theorem]{Hypothesis}

\theoremstyle{remark}
\newtheorem{remark}[theorem]{Remark}

\newtheorem{notation}[theorem]{Notation}
\newtheorem{context}[theorem]{Context}

\newcommand{\LFP}{{\rm LFP}}

\newcommand{\dn}{{\rm dn}}

\newcommand{\graph}{{\rm graph}}

\newcommand{\suc}{{\rm suc}}

\newcommand{\ind}{{\rm ind}}

\newcommand{\Prob}{{\rm Prob}}
\newcommand{\lev}{{\rm lev}}

\newcommand{\rang}{{\rm rang}}
\newcommand{\height}{{\rm height}}

\newcommand{\dg}{{\rm dg}}
\newcommand{\gr}{{\rm gr}}

\newcommand{\length}{{\rm length}}
\newcommand{\ggr}{{\rm gr}}

\newcommand{\rest}{{\restriction}}
\newcommand{\dom}{{\rm dom}}

\newcommand{\then}{{\underline{then}}}
\newcommand{\when}{{\underline{when}}}

\newcommand{\Iff}{{\underline{iff}}}
\newcommand{\mn}{{\medskip\noindent}}
\newcommand{\sn}{{\smallskip\noindent}}

\newcommand{\cE}{{\mathscr E}}

\newcommand{\cL}{{\mathscr L}}

\newcommand{\cF}{{\mathscr F}}

\newcommand{\bbL}{{\mathbb L}}
\newcommand{\bbN}{{\mathbb N}}

\newcommand{\bbQ}{{\mathbb Q}}
\newcommand{\bbR}{{\mathbb R}}

\newcommand{\varp}{{\varepsilon}}

\newcount\skewfactor
\def\mathunderaccent#1#2 {\let\theaccent#1\skewfactor#2
\mathpalette\putaccentunder}
\def\putaccentunder#1#2{\oalign{$#1#2$\crcr\hidewidth
\vbox to.2ex{\hbox{$#1\skew\skewfactor\theaccent{}$}\vss}\hidewidth}}

\newenvironment{PROOF}[2][\proofname.]
   {\begin{proof}[#1]\renewcommand{\qedsymbol}{\textsquare$_{\rm #2}$}}
   {\end{proof}}

\begin{document}

\title [Failure of 0-1 law for random graph in strong logics]
       {Failure of 0-1 law for sparse random graph \\
in strong logics}

\dedicatory  {Dedicated to Yuri Gurevich on the Occasion of his 75th
Birthday}

\author {Saharon Shelah}
\address{Einstein Institute of Mathematics\\
Edmond J. Safra Campus, Givat Ram\\
The Hebrew University of Jerusalem\\
Jerusalem, 9190401, Israel\\
 and \\
 Department of Mathematics\\
 Hill Center - Busch Campus \\ 
 Rutgers, The State University of New Jersey \\
 110 Frelinghuysen Road \\
 Piscataway, NJ 08854-8019 USA}
\email{shelah@math.huji.ac.il}
\urladdr{http://shelah.logic.at}
\thanks{This work was
partially supported by European Research Council grant 338821.
\
Publication 1062 on Shelah's list.
\
The author thanks Alice Leonhardt for the beautiful typing.
First typed Aug. 28, 2012.}






\date{March 24, 2017}

\begin{abstract}
Let $\alpha \in (0,1)_{\bbR}$ be irrational and $G_n =
G_{n,1/n^\alpha}$ be the random graph on $[n]$ with edge probability
$1/n^\alpha$; we know that it satisfies the 0-1 law for first order logic.
We deal with the failure of the 0-1 law for 
stronger logics: $\bbL_{\infty,\bold k},\bold k$ large enough
and the LFP, least fix point logic.
\end{abstract}

\maketitle
\numberwithin{equation}{section}
\setcounter{section}{-1}
\newpage

\section {Introduction}
\bigskip

\subsection {The Question} \
\bigskip

Let $G_{n,p}$ be the random graph with set of nodes $[n] =
\{1,\dotsc,n\}$, each edge of probability $p \in [0,1]_{\bbR}$, the
edges being drawn independently, see $\boxplus_2$ below.  On 0-1 laws 
(and random graphs)
see the book of Spencer \cite{Spe01} or Alon-Spencer \cite{AlSp08},
in particular on the behaviour of the random graph $G_{n,1/n^\alpha}$ for
$\alpha \in (0,1)_{\bbR}$ irrational.
On finite model theory see Flum-Ebbinghaus \cite{EbFl06}, e.g. on the
logic $\bbL_{\infty,\bold k}$ (see \S1) and on LFP 
(least fixed point\footnote{There are some variants, but
  those are immaterial for our perspective}) logic.  A characteristic example
of what can be expressed by it is ``in the graph $G$ there is a path
from the node $x$ to node $y$"; this is closed to what we use.
We know that $G_{n,p}$, i.e. $p$ constant satisfies the 0-1 law for 
first order logic (proved independently by Fagin
\cite{Fa76} and Glebskii-et-al \cite{GKLT}).  This holds also for many
stronger logics like $\bbL_{\infty,\bold k}$ and LFP logic.  If
$\alpha \in (0,1)_{\bbR}$ is irrational, the 0-1 law holds for
$G_{n,(1/n^\alpha)}$ and first order logic, see e.g. \cite{AlSp08}. 

The question we address is whether this holds also for stronger logics as
above.  Though our main aim is to address the problem for the case of graphs, the proof seems
more transparent when we have two random graph relations (we make them
directed graphs just for extra transparency).  So here we shall deal with two
cases A and B.  In Case A, the usual graph, we have to show that there
are (just first order) formulas $\varphi_\ell(x,y)$ for $\ell=1,2$
with some special properties (actually also $\varphi_0$), 
see Claim \ref{a2}.  For Case B, those
formulas are $R_\ell(x,y),\ell=1,2$, the two directed graph relations.  
Note that (for Case B), the
satisfaction of the cases of the $R_\ell$ are decided directly by the
drawing and so are independent, whereas for Case A there are (small)
dependencies for different pairs, so the probability estimates are
more complicated.

In the case of constant probability $p \in (0,1)_{\bbR}$,  
the 0-1 law is strong: it is obtained by proving elimination
of quantifier and it works also for stronger logics:
$\bbL_{\infty,\bold k}$ (see \S2) and so also for the LFP logic
$\bbL_{\LFP}$.  Another worthwhile case is:
\mn
\begin{enumerate}
\item[$\boxplus_1$]  $G_{n,1/n^\alpha}$ where $\alpha \in (0,1)_{\bbR}$;
  so $p_n = 1/n^\alpha$.
\end{enumerate}
\mn
Again the edges are drawn independently but the probability depends on
$n$.

The 0-1 law holds if $\alpha$ is irrational, but we have elimination
of quantifiers only up to (Boolean combinations of) existential formulas.  Do we have 0-1 law
also for those stronger logics?  We shall show that by proving
that for some so called scheme of interpretation $\bar\varphi$, for
random enough $G_n,\bar\varphi$ interpret number theory up to $m_n$
where $m_n$ is not too small, e.g. $m_n \ge \log_2 \log_2(n)$.

A somewhat related problem asks whether for some logic the 0-1 law
holds for $G_{n,p}$ (e.g. $p=\frac 12$) \underline{but} does not have
the elimination of quantifiers, see \cite{Sh:1077}.

We now try to \underline{informally} describe the proof,
 naturally concentrating on case B. 

Fix reals $\alpha_1 < \alpha_2$ from $(0,\frac 14)_{\bbR}$ for transparency,  
so $\bar\alpha = (\alpha_1,\alpha_2)$ letting $\alpha(\ell) = \alpha_\ell$; 
\mn
\begin{enumerate}
\item[$\boxplus_2$]  let the random digraph 
$G_{n,\bar\alpha} = ([n],R_1,R_2) = ([n],R^{G_{n,\bar\alpha}}_1,
R^{G_{n,\bar\alpha}}_2)$ with $R_1,R_2$
irreflexive 2-place relations drawn as follows:
\sn
\begin{enumerate}
\item[$(a)$]  for each $a \ne b$, we draw a truth value for
$R_2(a,b)$ with probability $\frac{1}{n^{1-\alpha_2}}$ for yes
\sn
\item[$(b)$]  for each $a \ne b$, we draw a truth value for
$R_1(a,b)$ with probability $\frac{1}{n^{1+\alpha_1}}$ for yes
\sn 
\item[$(c)$]  those drawings are independent.
\end{enumerate}
\end{enumerate}
\mn
Now for random enough digraph $G= G_n= G_{n,\bar\alpha} = 
([n],R_1,R_2)$ and node $a \in G$; we try to  define the set
$S_k = S_{G,a,k}$ of nodes of $G$ not from $\cup\{S_m:m< k\}$ 
by induction on $k$ as follows:

For $k=0$ let $S_k = \{a\}$.
Assume $S_0,\dots,S_k$ has been chosen, and we shall choose $S_{k+1}$.
\mn
\begin{enumerate}
\item[$\boxplus_3$]   For $\iota =1,2$ 
 we ask: is there an $R_\iota$-edge $(a,b)$ 
with $a \in S_k$ and $b$ not from $\cup\{S_m:m \le k\}$?
\end{enumerate}
\mn
If the answer is no for both $\iota =1,2$
we stop and let ${\height}(a, G)= k$.
If the answer is yes for $\iota= 1$, we let $S_{k+1}$
be the set of $b$ such that for some $a$
the pair $(a,b)$ is as above for $\iota =1$.
If the answer is no for $\iota =1$ but yes
for $\iota=2$ we define $S_{k+1}$ similarly
using $\iota=2$.

Let the height of $G$ be $\max\{\height(a,G):a \in G\}$.

Now we can prove that for every random enough 
$G_n$, for $a \in G_n$ or easier - for most
$a \in G_n$, for every not too large $k$ we have:
\mn
\begin{enumerate}
\item[$\boxplus_4$]  $S_{G_n,a,k}$ is on the one hand not empty
and on the other hand with $\le n^{2 \alpha_2}$ members.
\end{enumerate}
\mn
This is proved by drawing the edges not all at once
but in $k$ stages. In stage $m \le k$ we already can compute 
$S_{G_n,a,0},\dots,S_{G_n,a,m}$ and we
have already drawn all the $R_1$-edges
and $R_2$-edges having first node in $S_{G_n,a,0} \cup  \dots \cup 
S_{G_n,a,m-1}$; that is for every such pair 
$(a,b)$ we draw the truth values of $R_1(a,b),R_2(a,b)$. 
For $m=0$ this is clear.  So arriving to $m$ we can 
draw the edges having the first node in $S_m$
and not dealt with earlier, and hence can compute $S_{m+1}$.

The point is that in the question $\boxplus_3$ above, if the 
answer is yes for $\iota = 1$ 
then the number of nodes in $S_{m+1}$ will be small, almost surely
smaller than in $S_m$ because its expected value is $|S_m| \cdot |[n]
- \bigcup\limits_{\ell \le m} S_\ell| \cdot
\frac{1}{n^{1+\alpha_1}} \le n^{1+2 \alpha_2-(1+\alpha_1)} = n^{2
  \alpha_2 - \alpha_1}$ and the drawings are independent so except for
an event of very small probability this is what will occur.  
Further, if for $\iota =1$ 
the answer is no but for $\iota =2$ the answer is yes then 
almost surely $S_m$ is smaller than a number near $n^{\alpha_1}$ 
but it is known that the $R_2$-valency of any node of $G_n$ is
near to $n^{\alpha_2}$.  Of course, the ``almost surely" is such that the
probability that at least one undesirable event mentioned above occurs
is negligible.

So the desired inequality holds. 

By a similar argument, if we stop at $k$ then there is no $R_2$-edge
from $S_k$ into $[n] \backslash (S_0 \cup \ldots S_k)$ so the expected
value is $\ge |S_k| \cdot (n - \sum\limits_{\ell \le k} (S_k)) \cdot
\frac{1}{n^{1-\alpha_2}}$ hence in $S_0 \cup \dots \cup S_k$ 
there are many nodes, e.g. at least near $n/2$
by a crude argument.
As each $S_m$ is not too large necessarily
the height of $G_n$ is large.

The next step is to express in our logic the relation 
$\{(a_1,b_1,a_2,b_2)$: for some $k_1,k_2$  we have 
$b_1 \in S_{G_n,a_1,k_1},b_2 \in S_{G_n,a_2,k_2}$ and $k_1 \le k_2\}$.

By this we can interpret a linear order
with $\height(G_n)$ members.
Again using the relevant logic this suffice 
to interpret number theory up to this
height. Working more we can define a linear order 
with $n$ elements, so can essentially
find a formula ``saying" $n$ is even (or odd).

For random graphs we have to work harder: instead
of having two relations we have two formulas;
one of the complications is that
their satisfaction for the relevant pairs are not fully independent.

In \cite{Sh:1096} we shall deal with the strong failure of the 0-1 for Case
A, i.e. $G_{n,p^\alpha}$, (e.g. can ``express" $n$ is even) and also
intend to deal with the $\alpha$ rational case.
The irrationality can be replaced by discarding few exceptions.

We thank the referee for helping to improve the presentation.
\bigskip
 
\subsection {History} \
\bigskip

The history is non-trivial having non-trivial opaque points.  I 
have a clear memory of the events but vague on the exact statements and
more so on the proof
(and a concise entry in my (private F-list, \cite{Sh:F159}))
that in January 1996, in a Conference in DIMACS, Monica McArthur 
gave a lecture claiming that the graph
$G_{n,\alpha}$ satisfies the 0-1 law not only for first order logic
(by Shelah-Spencer \cite{ShSp:304}) but also for a stronger logic.
Joel Spencer said this can be contradicted 
in a way he outlined.  I thought on this and saw 
further things (explain to her saying more and)
wrote them in a letter to Monica and Joel.  I understood that 
it was agreed that Monica would write a paper with us 
saying more but eventually she left academia.

As the referee found out, MacArthur's claim in \cite{McA97}
(DIMACS): failure of the law in
$\bbL^\omega_{\infty,\omega}$, but refers the proof to a paper in
preparation with Spencer that never appeared.  She claims also that
there is 0-1 law for $\bbL^k_{\infty,\omega}$ if $k = [1/\alpha]$,
referring again to the paper in preparation.  The later claim is not
contradicted by the results of this paper.  Lynch \cite{Lyn97}, refers also
to a joint paper with McArthur and Spencer that never appeared proving
that for the TC (= transitive closure) logic satisfies the 0-1 law.

Having sent Joel (in 2011) an earlier version of this paper, his
recollection of talking to Monica was 
that ``we hadn't really gotten a handle on the situation".

Discussing with Simi Haber, (December 2011) this question
arised again, trying to recollect it was 
not clear to me what was the logic (inductive
logic? $\bbL_{\infty,k}$?)  Looking at it again, I saw a proof
for the logic $\bbL_{\infty,k}$. 
No trace of the letter or the notes mentioned above were found.  The only 
tangible evidence, is in an entry  
\cite{Sh:F159} from my F-list. 
Joel declined a suggestion that Haber, he and I 
will deal with it, and eventually also Haber left.

The notes on \S1 are from January 2012; for \S2 from Sept. 4, 2012;
revised in Nov/Dec. 2014 and expanded March 2015, June, 2015.

The intention was that it would appear in the Yurifest, commemorating
Yuri Gurevich's 75th birthday, but it was not in a final version in
time, so only a short version (with the abstract and \S(0A)) appear in
the Yurifest volume, \cite{Sh:1061}.
\bigskip

\subsection {Preliminaries} \
\bigskip

\begin{notation}  
\label{z2}
1) $n \in \bbN \backslash \{0\}$ will be used for $``G_n \in K_n$ random
  enough".

\noindent
2) $G,H$ denote graphs and $M,N$ denote more general structures = models. 

\noindent
3) $a,b,c,d,e$ denote nodes of graphs or elements of structures.

\noindent
4) $m,k,\ell$ denote natural numbers.

\noindent
5) $\tau$ denotes a vocabulary, $M$ a model with vocabulary
$\tau=\tau_M$ (see \ref{z2}(9),(10) below).

\noindent
6) $\cL$ denotes a logic, $\bbL$ is first order logic, 
so $\bbL(\tau)$ is first order language (= set of formulas) 
for the vocabulary $\tau$.
$\cL(\tau)$ is the language for the logic $\cL$ and the vocabulary $\tau$. 

\noindent
7) $\bbL_{\LFP}$ is the least fix point logic, short hand $\LFP$.

\noindent
8) 
\mn
\begin{enumerate} 
\item[$(a)$]  Let ${}^k A$ be the set of sequences $\eta$ of length $k$ of
members of $A$, i.e. $\eta = \langle
a_0,\dotsc,a_{k-1}\rangle$ where $\bigwedge\limits_{\ell < k} a_\ell \in A$,
so $a_\ell = \eta(\ell)$.
\sn
\item[$(b)$]  For a set $u$, e.g. of natural number let $\bar x_{[u]}
  = \langle x_s:s \in u \rangle$,
\sn
\item[$(c)$]  if $\varphi(\bar x_m,\bar y) \in \cL(\tau)$ and $M$ is a
$\tau$-model and $\bar b \in {}^{\ell g(\bar y)}M$ and $\bar x_m =
\langle x_i:i < m\rangle$ then $\varphi(M,\bar b) = \{\bar a \in {}^m M:M
\models \varphi[\bar a,\bar b]\}$.
\end{enumerate}
\mn
9) Let $\tau_{\gr}$ denote the vocabulary of graphs, but we may write
$\bbL(\graph)$ or $\cL(\graph)$ instead of
$\bbL(\tau_{\gr}),\cL(\tau_{\gr})$.  So $\tau_{\gr}$ consists of one
two-place predicate $R$, (below always interpreted as a symmetric
irreflexive relation).

\noindent
10) Let $\tau_{\dg}$ denote the vocabularies of bound directed graphs,
so it consists of two two-place predicates, below always
interpreted as irreflexive relations.  Let $\tau_{\bbN}$ be the
vocabulary of number theory, see \ref{z3}(1).

\noindent
11) We define the function $\log_*$ from $\bbR_{\ge 0}$ to $\bbN$ by: 

$\log_*(x)$ is 0 if $x<2$

$\log_*(x)$ is $\log_*(\log_2(x)) +1$ if $x \ge 2$

\noindent
12)  $|u|$ is the cardinality = the number of elements of a set $u$.
\end{notation}

\begin{explanation}
\label{z3}
1) Above recall that the vocabulary of the model = the structure 
$\bbN$, number theory, is the
set of symbols $\{0,1,+,\times,<\}$ where 0,1 are individual constants
(interpreted in $\bbN$ as the corresponding elements) and $+,\times$ are
two-place function symbols interpreted as $+^{\bbN},\times^{\bbN}$ the
two-place functions of addition and multiplication, and $<$ is a
two-place predicate (= relation symbol) interpreted as $<^{\bbN}$, the
usual order on $\bbN$.

\noindent
2) In general
\mn
\begin{enumerate}
\item[$(A)$]  a vocabulary is a set of predicates, individual constants
  and function symbols each with a given arity = number of places;
individual constants (like 0,1 above) are considered as 0-place function symbols
\sn
\item[$(B)$]  a $\tau$-model or a $\tau$-structure $M$ consists of:
\sn
\begin{enumerate}
\item[$(a)$]  its universe, $|M|$, a non-empty set of elements so
  $\|M\|$ is their number
\sn
\item[$(b)$]  if $P \in \tau$ is an $n$-place predicate, $P^M$ is a
  set of $n$-tuples of members of $M$
\sn
\item[$(c)$]  if $F \in T$ is an $n$-place function symbol then $F^M$
  is an $n$-place function from from $|M|$ to $|M|$.
\end{enumerate}
\end{enumerate}
\end{explanation}

\begin{definition}
\label{z3d}
Let $\tau$ be a finite vocabulary, for simplicity with predicates
only or we just consider a function as a relation.  
Here we use $\tau_{\ggr},\tau_{\dg}$ only except when we interpret.

\noindent
1) We say $\bar\varphi$ is in a $(\tau_*,\tau)$-scheme of interpretaion
\when \,: (if $\tau$ is clear from the context we may write
$\tau_*$-scheme)
\mn
\begin{enumerate}
\item[$(a)$]  $\bar\varphi = \langle \varphi_R(\bar x_{n_\tau(R)}):R
  \in \tau_* \cup \{=\}\rangle$ where $n_\tau(R)$ is the arity (=
  number of places) of $R$
\sn
\item[$(b)$]  $\varphi_R \in \bbL(\tau)$
\sn
\item[$(c)$]  $\varphi_=(x_0,x_1)$ is always an equivalence relation
  on $\{y:(\varphi(y,y)\}$; if $\varphi_=$ is $(x_0=x_1)$ then we may omit it.
\end{enumerate}
\mn
2) For a $\tau$-model $M$ (here a graph or diagram) and $\bar\varphi$
as above, let $N=N_{M,\bar\varphi}$ be the following structure:
\mn
\begin{enumerate}
\item[$(a)$]  $|N|$ the set of elements of $N$, is $\{a/\varphi_=(M):a \in
  M$ and $M \models \varphi_=(a,a)\}$; note that $\varphi(M)$ is an
  equivalence relation on $\{a:M \models \varphi_=(a,a)\}$
\sn
\item[$(b)$]  if $R \in \tau$ has arity $m$ then $R^N$, the
  interpretation of $r$ is $\{\langle a_\ell/\varphi_=(M):\ell < m:
M \models \bigwedge\limits_{\ell < m} \varphi_=(a_\ell,a_\ell)
  \wedge \varphi_R(a_0,\dotsc,a_{m-1}\}$ so $a_0,\dotsc,a_{m-1} \in M\}$.
\end{enumerate}
\end{definition}

\noindent
Recall that here ``for every random enough $G_n$" is a central notion.
\begin{definition}
\label{z3f}
1) A 0-1 context consists of:
\mn
\begin{enumerate}
\item[(a)]  a vocabulary $\tau$, here just the one of graphs or double
  directed graphs, see \ref{z2}(5),(9),(10)
\sn
\item[(b)]  for each $n,K_n$ is a set of $\tau$-models with set of
  elements = nodes $[n]$, in our case graphs or double driected
  graphs
\sn
\item[(c)]  a distribution $\mu_n$ on $K_n$, i.e. $\mu_n:K_n
  \rightarrow [0,1]_{\bbR}$ satisfying $\Sigma\{\mu_n(G):G \in
  K_n\}=1$
\sn
\item[(d)]  the random structure is called $G_n = G_{\mu_n}$ and we
  tend to speak on $G_{\mu_n}$ rather than on the context.
\end{enumerate}
\mn
2) For a given 0-1 context, let ``for every random enough $G_n$ we
have $G_n \models \psi$, i.e. $G$ satisfies $\psi$" and ``if $G_n$ is
random enough then $\psi$", etc. means that the sequence $\langle \Prob(G_n
\models \psi):n \in \bbN\rangle$ converge to 1; of course, $\Prob(G_n
\models \psi) = \Sigma\{\mu_n(G):G \in K_n$ and $G \models \psi\}$.

\noindent
3) For $\bar p = \langle p_n =p(n):n\rangle$ a sequence of
probabilities, $G_{n,\bar p}$ is the case $K_n =$ graphs on $[n]$ and
we draw the edges independently
\mn
\begin{enumerate}
\item[(a)]  with probability $p$ when $\bar p$ is constantly $p$,
  e.g. $\frac 12$, and
\sn
\item[(b)]  with probability $p(n)$ or $p_n$ when $p$ is a function
  from $\bbN$ to $[0,1]_{\bbR}$.
\end{enumerate}
\end{definition}

\noindent
Below we add the second context because for it the proof is more
transparent.
\begin{context}
\label{z4}
1) \underline{Case A}:
\mn
\begin{enumerate}
\item[(a)]  $a \in (0,1)$ is irrational
\sn
\item[(b)]  $p_n = 1/n^\alpha$.
\end{enumerate}
\mn
2) \underline{Case B}: 
\smallskip

\noindent
${\bar\alpha}^*  = (\alpha^*_1,\alpha^*_2)$ where
$\alpha_1^*,\alpha_2^*  \in (0,1/4)$ are irrational numbers, (natural
to add linearly independent over $\bbQ$) hence
$0 < \alpha^*_1 < \alpha^*_2 < \alpha^*_2 + \alpha^*_2  <  1/2$ 
and let $\alpha^*_0 = \alpha^*_1$. 
\end{context}

\begin{definition}  
\label{z8}
For Case A:

\noindent
1) Let $K^1 := \bigcup\limits_{n} K^1_n$ where we let
$K^1_n$ be the set of graphs $G$ on $[n] = \{1,\dotsc,n\}$ 
so $R^G \subseteq \{\{i,j\}:i \ne j \in [n]\}$.

\noindent
2) For $\alpha \in (0,1)_{\bbR}$ let 
$G_n = G_{n;\alpha}$ be the random graph on $[n]$ with the
   probability of an edge being $1/n^\alpha$ and the drawing of the
   edges being independent.

\noindent
3) Let $\mu_n = \mu_{n;\alpha}$ be the corresponding distribution on
 $K^1_n$; so $\mu_n:K^1_n \rightarrow [0,1]_{\bbR}$ and $1 = 
\Sigma\{\mu_n(M):M \in K_n\}$, in fact, $\mu_n(G) =
   (1/n^\alpha)^{|R^G|} \times (1-1/n^\alpha)^{\binom n2-|R^G|}$.
\end{definition}

\begin{convention}
\label{z9}
Writing $K_n$ means we intend $K^1_n$ or $K^2_n$ (see below),
similarly $G_n$ is $G_{n,\alpha}$ if Case A and $G_{n,\bar\alpha}$ if
Case B and similarly $K$ is $K^1$ or $K^2$.
\end{convention}

\noindent
The more transparent related case is the following
\begin{definition}
\label{z13}
On Case B, for $G_{n;\bar\alpha}$:

\noindent
1) Recall $\tau_{\dg}$ is the 
vocabulary $\{R_1,R_2\}$ intended to be two directed
graph relations.

\noindent
2) Let $K^2 = \bigcup\limits_{n} K^2_n$ where we let
$K^2_n = \{G:G = ([n],R^G_1,R^G_2)$ satisfying 
$([n],R^G_\ell)$ is a directed graph for $\ell=1,2$; 
we may write $R_\ell$ instead of
$R^G_\ell$ when $G$ is clear from the context$\}$.  
We assume\footnote{We may change the definition of $K^2_n$ by
  requiring $R^G_1 \cap R^G_2 = \emptyset$, this makes little
  difference.  We could further demand $R_\ell$ is asymmetric,
  i.e. $(a,b) \in R^G_\ell \Rightarrow (b,a) \notin R^G_0$, again this
  makes little difference.} irreflexivity, i.e. $(a,a) \notin R^G_\ell$ 
but allow $(a,b),(b,a) \in R^G_\ell$.

\noindent
3) For reals $\alpha_1 < \alpha_2$ from $(0,\frac 14)_{\bbR}$, say from
\ref{z4}(2) so $\bar\alpha = (\alpha_1,\alpha_2)$
 let $\alpha(\ell) = \alpha_\ell$; let the random model
 $G_{n;\bar\alpha} = ([n],R_1,R_2) = ([n],
R^{G_{n;\bar\alpha}}_1,R^{G_{n;\bar\alpha}}_2)$ with $R_1,R_2$ 
irreflexive relations be drawn as follows:
\mn
\begin{enumerate}
\item[$(a)$]  for each $a \ne b$, we draw a truth value for
  $R_2(a,b)$ with probability $\frac{1}{n^{1-\alpha_2}}$ for yes
\sn
\item[$(b)$]  for each $a \ne b$, we draw a truth value for
$R_1(a,b)$ with probability $\frac{1}{n^{1+\alpha_1}}$ for yes
\sn
\item[$(c)$]  those drawings are independent.
\end{enumerate}
\mn
4) We define the distribution $\mu_{n;\bar\alpha}$ as follows:
\mn
\begin{enumerate}
\item[$(a)$]  $\mu_n = \mu_{n;\bar\alpha} = 
\mu_{n;\alpha_1,\alpha_2}$ is the following
distributions on $K^2_n$:
\sn
\begin{enumerate}
\item[$\bullet$]  $\mu_n(G) = \mu^2_{n;\alpha_2}([n],R^G_1)) \cdot
 \mu^2_{n;-\alpha_1}([n],R^G_2))$ where
\sn
\item[$\bullet$]  $\mu^2_{n;\alpha}([n],R) =
  (\frac{1}{n^{1-\alpha}})^{|R^G|} \cdot
  (1-\frac{1}{n^{1-\alpha}})^{n(n-1)-|R^G|}$
\end{enumerate}
\sn
\item[$(b)$]  $G_n = G_{n;\bar\alpha} = G_{n;\alpha_1,\alpha_2}$
 denote a random enough $G \in K^2_n$ for $\mu_{n;\bar\alpha}$ 
so $n$ is large enough.
\end{enumerate}
\end{definition}

\begin{observation}
\label{z17}
For random enough (recalling \ref{z3f}(2))
$G_n = G_{n;\bar\alpha} = G_{n;\alpha_1,\alpha_2}$:
\mn
\begin{enumerate}
\item[$(a)$]  For $a \in [n]$, the expected value of the $R_2$-valency
  of $a$, that is, $|\{b:a R^G_2 b\}|$
 is $(n-1) \cdot \frac{1}{n^{1-\alpha(2)}} \sim n^{\alpha(2)}$;
\sn
\item[$(b)$]  for every random enough $G_{n;\bar\alpha}$
for every $a \in [n]$ this number is close 
enough to $n^{\alpha(1)}$, e.g.
\sn
\begin{enumerate}
\item[$\bullet_2$]  for some $\varp \in (0,\alpha_1)_{\bbR}$, the
  probability of the difference being $\ge n^{\alpha(1)(1-\varp)}$ for
  at least one $a \in [n]$, goes to zero with $n$
\end{enumerate}  
\sn
\item[$(c)$]  the expected number of $R_1$-edges
 is $n(n-1)/n^{+(1 + \alpha_1)} \sim n^{1-\alpha_1}$ hence the
 expected value of $|\{a:aR_1 b$ for some $b\}|$ is close to it
\sn
\item[$(d)$]  for every random enough $G_{n,\bar\alpha}$
the two numbers in (c) are close enough to
  $n^{1-\alpha_1}$ (similarly to (b)).
\end{enumerate}
\end{observation}

\begin{remark}
1) For $K^2$, this is a parallel of Claim \ref{a2} for $K^1$.

\noindent
2) Note that the Clause (a) does not imply clause (b) in \ref{z17}
because a priori the variance may be too large.
\end{remark}
\newpage

\section {On the logic $\bbL_{\infty,\bold k}$}

As the proof for $\bbL_{\infty,\bold k}$ is simpler and more
transparent than for $\LFP$, we shall explain it.

First, we try to define and then explain the logic
$\bbL_{\infty,\bold k}$ for $\bold k$ a finite number.  

For a vocabulary $\tau$, we
define the set $\bbL_{\infty,\bold k}(\tau)$ of formulas as the closure of
the set of atomic formulas under some operation similarly to first
order logic, but:
\mn
\begin{enumerate}
\item[$\bullet$]  we restrict ourselves to formulas having $< \bold k$
  free variables
\sn
\item[$\bullet$]  we allow arbitrary conjunctions and disjunctions
  (that is even infinite\footnote{As we consider only finite models,
    countable conjunctions and injunctions are enough.} ones)
\sn
\item[$\bullet$]  as in first order logic we allow negation
  $\neg\varphi$ and existential quantifier (on one variable) $\exists
  x \varphi(x,\bar y)$.
\end{enumerate}
\mn
So any formula in $\bbL_{\infty,\bold k}$ not just have only $< \bold k$
free variables but also every subformula has.

It may be helpful to recall the standard game which express
equivalence.  Recall (\ref{z3d}(1)) that for 
transparency we assume the vocabulary below has only predicates and is finite.
\mn
\begin{enumerate}
\item[$\boxplus$]  we say $\cF$ is an
  $(M_1,M_2)-\bbL_{\infty,\bold k}$-equivalence witness \when \, for some
  vocabulary $\tau$ with predicates only
\sn
\begin{enumerate}
\item[$(a)$]  $M_1,M_2$ are $\tau$-models
\sn
\item[$(b)$]  $\cF$ is a non-empty set of partial isomorphisms from $M_1$ to
  $M_2$
\sn
\item[$(c)$]  if $f \in \cF$ then $|\dom(f)| < \bold k$
\sn
\item[$(d)$]  if $f \in \cF,A \subseteq \dom(f),|A|+1 < \bold k,\iota
\in \{1,2\}$ and $a_\iota \in M_\iota$ \then\, there is $g$ such that
$g \in \cF,f \rest A \subseteq g$ and $\iota =1 \Rightarrow a_\iota \in
\dom(g)$ and $\iota = 2 \Rightarrow a_\iota \in \rang(g)$.
\end{enumerate}
\end{enumerate}
\mn
Now
\mn
\begin{enumerate}
\item[$\oplus_1$]  for $M_1,M_2$ as in (a) of $\boxplus$ above, 
the following are equivalent:
\sn
\begin{enumerate}
\item[$(a)$]  $M_1,M_2$ are $\bbL_{\infty,\bold k}$-equivalent,
  i.e. for every sentence $\psi \in \bbL_{\infty,\bold k}(\tau)$,
  i.e. a formula with no free variables, $M_1 \models \psi
  \Leftrightarrow M_2 \models \psi$
\sn
\item[$(b)$]  there is an $(M_1,M_2)-\bbL_{\infty,\bold
    k}$-equivalence witness $\cF$, i.e. as in $\boxplus$.
\end{enumerate}
\end{enumerate}
\mn
Also
\mn
\begin{enumerate}
\item[$\oplus_2$]   for $M_1,M_2,\cF$ as in $\boxplus$ above
  we have
\sn
\begin{enumerate}
\item[$(c)$]  if $k < \bold k,a_0,\dotsc,a_{k-1} \in M_1$ and $g \in
  \cF$ and $\{a_\ell:\ell < k\} \subseteq \dom(g)$ \then \, for
every formula $\varphi(x_0,\dotsc,x_{k-1}) \in 
\bbL_{\infty,\bold k}(\tau)$ we have
\[
M_1 \models \varphi[a_0,\dotsc,a_{k-1}] \Leftrightarrow M_2 \models
\varphi[g(a_0),\dotsc,g(a_{k-1})].
\]
\end{enumerate}
\end{enumerate}
\bigskip

\centerline {$* \qquad * \qquad *$}
\bigskip

Having explained the logic, how can we prove for it the failure of the
0-1 law?  Consider Case B where we have two kinds of edges, $R_1$ and
$R_2$.  Consider $\eta$ a sequence from ${}^k\{1,2\}$, see \ref{z2}(11) and
$a \ne b$.  There may be $(\eta,0,k)$-pre-paths from $a$ to $b$ in $G$, see
Definition \ref{a9}, i.e. $a = a_0,a_1,\ldots,a_k=b$ such that
$(a_\ell,a_{\ell +1})$ is an $R_{\eta(\ell)}$-edge for $\ell < k$.

Now depending on $\eta$ there may be many such pre-paths or few.  If
$\eta$ is constantly 2 and $k > 1/\alpha^*_2$ then there are many
such pre-paths - as fixing $a$ in $G_{n,\bar\alpha^*}$ the expected 
number of $b$'s for which there
is pre-$(\eta,0,k)$-paths from $a$ to $b$ is 1 for $k=0$, is $\approx
n^{\alpha^*_2}$ for $k=1$ is $\approx n^{2 \alpha^*_2}$ for $k=2$, etc., so
for $k > 1/\alpha^*_2$ it is every $b \in G_n$; not helpful.  If $\eta$ is
constantly 1, there are few such pre-paths and they are all short,
even $\le k$ for any random enough $G_n$, when $1 < \alpha^*_1 k$, not helpful.

But we may choose a ``Goldilock's'' $\eta$, that is, such that for
every initial segment of $\eta$ the expected number is not too large
and not too small.  This means that for some $a$ for every $k' \le k$ 
for some $b$ there is
such a pre-path but not too many.  
We need more so that we can define by a formula 
from $\bbL_{\infty,\bold k}$ the set
$S_{G_n,a,k'} := \{b$: there is such 
pre-path from $a$ to $b$ of length $k'$ but
not a shorter pre-path$\}$ and it is $\ne \emptyset$; 
moreover we can define the natural order on the set $\{S_{G_n,a,k}:k \le n\}$.
Fact \ref{a4} below indicates what kind of $\eta$'s we need, and we
use it proving \ref{a12}; however in later sections, because we have to
estimate the probabilities, we shall use only a closely related definition.

\begin{hypothesis}
\label{a0}
1) Case A  of \ref{z4} holds or Case B there holds.

\noindent
2) 
\mn
\begin{enumerate}
\item[$(a)$]  for case A: $\alpha^*_\ell,\varphi_\ell(x,y),
n^*_\ell$ for $\ell=0,1,2$ will be as in Claim \ref{a2} below
\sn
\item[$(b)$]  for case B: $\alpha^*_1,\alpha^*_2$ are as in \ref{z4}
  and $\varphi_\ell(x,y) = R_\ell(x,y)$ for $\ell=1,2$ and 
let $\alpha^*_0  = \alpha^*_1,\varphi_0(x,y) = \varphi_1(x,y)$
\sn
\item[$(c)$]  let $\bar\varphi = \langle \varphi_\ell(x,y):\ell =
  0,1,2\rangle$. 

\end{enumerate}
\end{hypothesis}

\begin{claim}  
\label{a2}
Assume we are in Case A. 
There are $\alpha^*_\ell,\varphi_\ell(x,y)$ and $\gamma^*_\ell$ for $\ell=0,1,2$ such that:
\mn
\begin{enumerate}
\item[$(a)$]  $0 < \alpha^*_1 < \alpha^*_0 < \alpha^*_2 < \alpha^*_2 +
  \alpha^*_2$ are reals with $\alpha^*_\ell \in (0,1/4)_{\bbR}$ and 
$\gamma^*_\ell \in \bbR_{>0}$
\sn
\item[$(b)$]  $\varphi_\ell(x,y)$ are first order formulas
  (in the vocabulary of graphs) even existential positive 
formulas such that $\varphi_\ell(x,y) \vdash x \ne y$ for random
enough $G_{n;\bar\alpha}$
\sn
\item[$(c)$]  if $G_{n;\bar\alpha}$ is random enough 
\then \, for every $a \in G_{n;\bar\alpha}$ the set 
$\varphi_2(G_{n;\bar\alpha},a)$ has $\approx \gamma^*_\ell
n^{\alpha^*_2}$ elements, i.e. for some 
$\varepsilon \in (0,1)_{\bbR}$, if $G_{n;\bar\alpha}$ is 
random enough, then for every $a \in [n]$, the number of
members of $\varphi_1(G_{n;\bar\alpha},a)$ belongs to the interval
  $(\gamma^*_\ell n^{\alpha^*_2} -
n^{\alpha^*_2(1-\varepsilon)},\gamma^*_2 n^{\alpha^*_2} +
  n^{\alpha^*_2(1 - \varepsilon)})$
\sn
\item[$(d)$]  if $\ell=0,1$ and 
$G_{n;\bar\alpha}$ is random enough \then \, $\{a \in [n]:\varphi_\ell
(G_{n;\bar\alpha},a) \ne \emptyset\}$ has 
$\approx \gamma^*_\ell n/n^{\alpha^*_\ell}$ members.
\end{enumerate}
\end{claim}

\begin{remark}
\label{a3}
We shall use not just the statements but 
also the proofs of \ref{a2}, \ref{a4}.
\end{remark}

\begin{PROOF}{\ref{a2}}  
Also here we shall use freely the analysis of $G_{n;\alpha}$ for
$\alpha \in (0,1)_{\bbR}$ irrational (see, e.g. \cite{AlSp08}).

Let $m^*_2,n^*_2$ be such that:
\mn
\begin{enumerate}
\item[$(a)$]   $n^*_2$ is large enough
\sn
\item[$(b)$]  $m^*_2 \le \binom{n^*_2 
}{2}$
\sn
\item[$(c)$]  $\alpha^*_2 := (n^*_2 -1) - \alpha m^*_2$ is positive
  but, e.g. $< \frac{1}{12}$.
\end{enumerate}
\mn
As $\alpha \in (0,1)_{\bbR}$ is irrational we can find such $m^*_2,n^*_2$.
Let $H^*_2$ be a random enough graph on $[n^*_2]$ with $m^*_2$ edges
such that $(1,2) \notin R^{H^*_2}$.  (Note that this ``random enough"
is just used for the existence proof).

We choose $n^*_1,m^*_1,H^*_1$ similarly except that
$- \alpha^*_1 := n^*_1 - 1 - \alpha m^*_1$ is negative with value close
enough, e.g. to $-\alpha^*_2/3$.  Lastly, we choose
$n^*_0,m^*_0,H^*_0$ similarly except that $- \alpha^*_0 = n^*_0-1-\alpha
m^*_6$ is negative and $\alpha^*_0 \in (\alpha^*_1,\alpha^*_2)$.

For $\ell=1,2$ let 
$\varphi_\ell(x,y) = (\exists \ldots x_i \ldots)_{i \in [n^*_\ell]}
(x=x_1 \wedge y = x_2 \wedge \bigwedge\{x_i Rx_j:i,j \in [m^*_\ell]$
satisfies $i R^{H^*_\ell} j\}) \wedge \bigwedge \{ x_i \not= x_j:
i \not= j \in [n^*_{\ell}]\}$.  

Now check clauses (a)-(d).  Clearly $\alpha^*_1,\alpha^*_2$
satisfy clause (a) and $\varphi_1,\varphi_2$ are as in clause (b).

For $\ell=1,2$, let $\gamma^*_\ell = 1$.  So for any $n$
large enough compared to $n^*_1,n^*_2$ and $a_1 \ne a_2 \in [n]$, the
set $\cF :=\{f:f$ is a one-to-one function from $[n^*_2]$ 
to $[n]$ such that $f(1) =a_1,f(2) = a_2\}$ has 
$\prod\limits_{i < n^*_2-2} (n-2-i) \sim n^{n^*_2-2}$ members.

For each $f \in \cF$ the probability of the event $\cE_f = ``f$ maps
every edge of $H^*_2$ to an edge of $G_{n,\alpha}$" is
$(\frac{1}{n^\alpha})^{m^*_2}$ so the expected value of $\{f \in
\cF:\cE_f$ occurs$\}$ is $\approx n^{n^*_2-\alpha m^*_2-2} =
\frac{1}{n^{1-\alpha^*_2}}$.  Clearly as in \ref{z17} the expected value is as
required in clause (c) and by the well known analysis of
$G_{n \alpha}$ (see, e.g. \cite{AlSp08}), clause (c) holds and see more in
\S4.

Clause (d) is proved similarly.
\end{PROOF}

\begin{fact}  
\label{a4}
There is a sequence $\eta \in {}^{\bbN}\{1,2\}$ such that: for every
$n > 0$,
\newline
$\gamma_n = |(\eta \rest n)^{-1}\{2\}| \alpha^*_2 - 
|(\eta \rest n)^{-1}\{1\}| \alpha^*_1$ belongs\footnote{Can and will use also
  other intervals and similar sequences} 
to $[\alpha^*_2 -\alpha^*_1,\alpha^*_2 + \alpha^*_2]_{\bbR}$.
\end{fact}

\begin{PROOF}{\ref{a4}}  
We choose $\eta(n)$ by induction on $n$.  Let $\eta(n)$ be 2 if
$\gamma_n \le \alpha^*_2$, e.g. $n=0$ and $\eta(n)$ be $1$ if $\gamma_n
> \alpha^*_2$.

Easily $\eta$ is as required.
\end{PROOF}

\begin{claim}  
\label{a7}
1) If $\eta$ is as in \ref{a4} \then \, for any $m$ and every random
 enough $G_n$, there is an $(\eta,m)$-path in $G_n$, see below.

\noindent
2) Moreover, also there is an $(\eta,\varepsilon \lfloor \log(n) 
\rfloor$-path and
even an $(\eta,\lfloor n^\varepsilon \rfloor)$-path for appropriate
$\varepsilon \in \bbR_{>0}$.
\end{claim}

\begin{PROOF}{\ref{a7}}
As in \cite{AlSp08} on $G_{n,1/n^\alpha}$ 
and see more in \S3.
\end{PROOF}

\begin{definition}  
\label{a9}
1) A sequence $\bar a = \langle a_\ell:\ell \in [m_1,m_2]\rangle$ 
of nodes, that is, of members of $G \in K_n$ is called a 
pre-$(\nu,m_1,m_2)$-path, if $m_1=0$ we may omit it, \when \,:
\mn
\begin{enumerate}
\item[$(a)$]  $\nu$ is a sequence of length $\ge m_2$ and $i < \ell
  g(\nu) \Rightarrow \nu(i) \in \{1,2\}$
\sn
\item[$(b)$]  if $\ell \in \{m_1,m_1+1,\dotsc,m_2-1\}$ then $G \models
\varphi_{\nu(\ell)}(a_\ell,a_{\ell+1})$.
\end{enumerate}
\mn
2) Above we say ``$(\nu,m_1,m_2)$-path" \when \, in addition:
\mn
\begin{enumerate}
\item[$(c)$]  if $m_1 \le \ell_1 < \ell_2 \le m_2$ and 
$\langle a'_\ell:\ell \in [m_1,\ell_2)\rangle$ is a 
pre-$(\nu,m_1,\ell_2)$-path \then \, $a'_{m_1} = a_{m_1} \wedge a'_{\ell_2} 
= a_{\ell_2} \Rightarrow a'_{\ell_1} = a_{\ell_1}$
\sn
\item[$(d)$]  if $m_1 \le \ell_1 < \ell_2 \le m_2$ \then \,
$a_{\ell_1} \ne a_{\ell_2}$.
\end{enumerate}
\mn
3) We say ``$\bar a$ is a (pre)-$(\nu,m_1,m_2)$-path from $a$ to $b$"
\when \, in addition $a_{m_1} = a \wedge a_{m_2} = b$.
\end{definition}

\begin{remark}
\label{a11}
1) Note that if $\langle a_\ell:\ell \le m\rangle$ is a
   pre-$(\nu,m)$-path, it is possible that $\ell_1 + 1 < \ell_2 \le m$
   and $a_{\ell_1} = a_{\ell_2}$.  For a $(\nu,m)$-path this is
   impossible.

\noindent
2) In \ref{a9}(2)(c), really the case $\ell_2 = m_2$ suffice.

\noindent
3) We use the ``$\log(n)$" Case in \ref{a7}(2), but having $\log(\log(n))$
   and even much less has no real affect on the proof.
\end{remark}

\begin{conclusion}  
\label{a12}
Let $\bold k \ge \max\{
n^*_0 , 
n^*_1,n^*_2\}$ \then \, $G_n$ fails the
0-1 law for $\bbL_{\infty,\bold k}$.
\end{conclusion}

\begin{remark}
\label{a11}
1) Note that if $\langle a_\ell:\ell \le m\rangle$ is a
   pre-$(\nu,m)$-path, it is possible that $\ell_1 + 1 < \ell_2 \le m$
   and $a_{\ell_1} = a_{\ell_2}$.  For a $(\nu,m)$-path this is
   impossible.

\noindent
2) In \ref{a9}(2)(c), really the case $\ell_2 = m_2$ suffice.

\noindent
3) We use the ``$\log(n)$" case in \ref{a7}(2), but having $\log \log(n)$
   and even much less has no real affect on the proof.

Note that we rely on \ref{a7}(2) but we prove more in \S3.
\end{remark}

\begin{PROOF}{\ref{a12}}  
For a finite graph $G$ and $\eta$ as in \ref{a4} or
any $\eta \in {}^{\bbN}\{1,2\}$ let $\length_\eta(G)$ be the maximal $m$ such
that there is an $(\eta,m)$-path in $G$.

Now consider the statement
\mn
\begin{enumerate}
\item[$\oplus$]  there is a sentence $\psi_m = \psi_{\eta,m} \in
\bbL_{\infty,\bold k}(\tau_{\ggr})$ such that for a finite 
graph $G,G \models \psi_m$ iff there is an $(\eta,m)$-path in $G$.
\end{enumerate}
\mn
Why $\oplus$ is enough?  Because then we let

\[
\psi = \bigvee\{(\psi_m \wedge \neg \psi_{m+1}):m \ge 10 \text{ and }
(\log_*(m) \text{ is even})\}
\]

\mn
where $\log_*(m)$ is essentially the inverse of the tower function,
see \ref{z2}(3).  Note that using \ref{a4}, \ref{a7}(2), 
of course, we should be able to say much more.

Why $\oplus$ is true?  First, we define the formula
$\psi_{m_1,m_2}(x,y)$ for $m_1 \le m_2$ by induction on $m_2-m_1$ as follows:
\mn
\begin{enumerate}
\item[$(*)_1$]  if $m_1 = m_2$ it is $x=y$
\sn
\item[$(*)_2$]  if $m_1 < m_2$ it is $(\exists
  x_1)[\varphi_{\eta(m_1)}(x 
,x_1) \wedge \psi_{m_1 +1,m_2}(x_1,y)]$.
\end{enumerate}
\mn
So
clearly 
\mn
\begin{enumerate}
\item[$(*)_3$]  if $G \in K_n$ and $a,b \in [n]$ then $G \models
  \psi_{m_1,m_2}(a,b)$ \Iff \, there is a pre-$(\eta,m_1,m_2)$-path
  from $a$ to $b$.
\end{enumerate}
\mn
Second, we define $\psi'_{m_2}(x,y)$ 
as 
 $\psi_{0,m_2}(x,y) \wedge
\bigwedge\limits_{\ell_1 < \ell_2 \le m_2} \neg(\exists
z'_1,z''_1,z_2)[z'_1 \ne z''_2 \wedge \psi_{0,\ell_1}(x,z'_1) \wedge
  \psi_{0,\ell_1}(x,z''_1) \wedge \psi_{\ell_1,\ell_2}(z'_1,z_2)
  \wedge \psi_{\ell_1,\ell_2}(z''_1,z_2) \wedge
  \psi_{\ell_2,m_2}(z_2,y)]$.

This just formalizes \ref{a9}(2)(c) so
\mn
\begin{enumerate}
\item[$(*)_4$]   $G \models \psi'_{m_2}(a,b)$ \Iff \, there is an
  $(\eta,m_2)$-path from $a$ to $b$.
\end{enumerate}
\mn
As said above this is enough.  Note that complicating the sentence we
may weaken the demand on $G_n$.
\end{PROOF}
\newpage

\section {The LFP Logic}

In this section we try to interpret an initial segment of number
theory in a random enough $G \in K_n$, i.e. with set of nodes $[n]$.
In Definition \ref{b4} for $G \in K_n$ and $a \in G$ we define a model
$M_{G,a}$.  Now in $M \in \bold M_{G,a_*}$, the equivalence
classes of $E^M$ represent natural numbers.  Concentrating on Case B,
starting with $\{a_*\}$ as the first equivalence class, usually its set of
$R_2$-neighbors will be the second equivalence class.  Generally, if
for an equivalence class $a/E^M$ we let the next one be the set
$\suc(a/E^M) = \{b \in G:R_2(a',b)$ for some $a' \in a/E^M\}$, 
then we expect that $\suc(A/E^M)$ has $\approx|a/E^M| 
\cdot n^{\alpha_2}$ members.  So if
we continue in this way, shortly we get the equivalence classes cover
essentially all the nodes of $G$.  Hence we try to sometimes use the
$R_1$-neighbors instead of the $R_2$-neighbors, but when?  For
$\bbL_{\infty,\bold k}(\tau_*)$ we can decide a priori, e.g.use $\eta$
as in \ref{a4} and the proof of \ref{a12} so that the
expected number will be small.  But for LFP logic this is not
clear, so we just say: use the $R_1$-neighbors if there is at least
one and the $R_2$-neighbors otherwise, so this is close to what is
done in \ref{a4}, \ref{a7}, \ref{a12} but not the same.

For case A we use $\varphi_\ell$-neighbors instead of $R_\ell$-neighbors
for $\ell=1,2$ except that the question on existence is for
$\varphi_0$-neighbors.

How do we from equivalence relations and the successor relation
re-construct the initial segment of number theory?  This is 
exactly the power of
definition by induction. 

Naturally we need just
\mn
\begin{enumerate}
\item[$\boxplus$]  letting $\height(G)$ be the maximal number of
  $E^M$-equivalence classes for $M \in \bold M_{G,a},a \in G$, we have:
\sn
\begin{enumerate}
\item[$(*)$]  for every $m$, for every random enough $G_n,m \le
  \height(G)$, moreover letting $f:\bbN \rightarrow \bbN$ be 
$f(n) = \log_*(n)$ for every random enough $G_n,f(n) \le
  \height(G_n)$. 
\end{enumerate}
\end{enumerate}
\mn
For failure of 0-1 laws, $\boxplus$ is enough, but we may wish to prove a
stronger version, say finding a sentence $\psi$ which for every random
enough $G_n,G_n$ satisfies $\psi$ \Iff \, $n$ is even.

We intend to return to it elsewhere; but for 
now note that for a set $A \subseteq G$ we can define 
($R$ is $R_2$ for Case B, $\varphi_2$ for Case A) 
$c \ell_{G_n}(A) = \{b:b \in A$ or $b \in G \backslash A$ but for no 
$c \in G \backslash A \backslash \{a\}$ do
we have $(\forall x \in A)[R(x,c) \equiv R(x,b)]\}$.  Now 
from a definition of a linear order on $A$ we can derive one
on $c \ell_{G_n}(A)$.  We can replace $R$ by any formula
$\varphi(x,y)$.  Now in our context, if we know that, with
parameters, we define such $A$ of size $\approx n^\varepsilon$ for appropriate
$\varepsilon$ then we can define a linear order on $[n]$; why there is
such $A$?  because if $M \in \bold M_{G,a_*}$ and $k$ is not too
large \then \, there is $b \in M,\lev(b,M)=k$ such that there is a unique
maximal $<^M_2$-path from $a_*$ to $b$.

For $\bbL_{\infty,\bold k}$ this is much easier.

\begin{context}
\label{b2}
(A) or (B):
\mn
\begin{enumerate}
\item[$(A)$]  (case A of \ref{z4}) the vocabulary $\tau_*$ is 
$\tau_{gr}$, the one for the class of graphs, 
$\varphi_\ell(x,y),\ell=0,1,2$ are as in \ref{a2} so $\in
\bbL(\tau_*)$ and $\bar\alpha^* = (\alpha^*_0,\alpha^*_1,\alpha^*_2)$, is
as there, $G = G_n = G_{n;\alpha},K_n = K^1_n$ as in Definition
\ref{a4}, 
\sn
\item[$(B)$]  (case B of \ref{z4}) $\tau_* = \tau_{\dg} =
\{R_1,R_2\},K_n = K^2_n$ and
$\bar\alpha^* = (\alpha^*_1,\alpha^*_2)$ 
are as in \ref{z4}(2) and $G_n = G_{n;\bar\alpha^*}$ 
and $\varphi_\ell(x,y) = x R_\ell y$ for
  $\ell=1,2$, with $G_n,K^2_n$ as in Definition \ref{z13} and let
  $\alpha^*_0 = \alpha^*_1,\varphi_0(x,y) = \varphi_1(x,y)$.
\end{enumerate}
\end{context}

\begin{definition}
\label{b4}
For $G \in K_n$ and $a_* \in G$ we define $\bold M = \bold M(G,a_*) 
= \bold M_{G,a_*}$ as the set of $\tau_1$-structures of $M$
such that (the vocabulary $\tau_1$ is defined implicitly):
\mn
\begin{enumerate}
\item[$(A)$]
\begin{enumerate}
\item[(a)]   the universe of $M$ is $P^M \subseteq [n]$
\sn
\item[(b)]  $c^M_* = a_*$, so $c_*$ is an individual
  constant from $\tau_ 1 $ 
\sn
\item[(c)]  $E^M$ is an equivalence relation on $M$
\sn
\item[(d)]  $<_1$ is a linear order on $P^M/E^M$, i.e.
\begin{enumerate}
\item[$(\alpha)$]  $a_1 E^M a_2 \wedge b_1 E^M b_2 \wedge a_1 <^M_1
  b_1 \Rightarrow a_2 <^M_1 b_2$
\sn
\item[$(\beta)$]  for every $a,b \in P^M$ 
exactly one of the following holds: 
$a <^M_1 b,b <^M_1 a$ and $a E^M b$
\end{enumerate}
\sn
\item[(e)]  $<^M_2$ is a partial order included in $<^M_1$
\sn
\item[(f)]
\begin{enumerate}
\item[$(\alpha)$]  $a_*/E^M$ is a singleton and $a_*$ is
  $<^M_2$-minimal, i.e. $b \in M \backslash \{a_*\} \Rightarrow a <^M_2 b$
\sn
\item[$(\beta)$]   if $a <^M_1 b <^M_1 c$ and $a <^M_2
  c$ then for some $b' \in b/E^M$ we have $a <^M_2 b' <^M_2 c$
\sn
\item[$(\gamma)$]   if $a,b \in M$ and $b/E^M$ is the
  immediate successor of $a/E^M$ (by $<^M_1$) \then \, for
 some $a' \in a/E^M$, we have $a' <^M_2 b$
\end{enumerate}
\sn
\item[(g)]  if $b \in M$ is a $<^M_2$-immediate
  successor of $a \in M$ (i.e. $a <^M_2 b$ 
 and $\neg(\exists y)(a <_2 y <_2 b)$, equivalently, 
$\neg(\exists y)(a <_1 y <_1 b))$
\then \, for some $\iota \in \{1,2\}$ 
we have $G \models \varphi_\iota[a,b]$
\sn
\item[(h)]  $P_0,P_1,P_2 = P_+,P_3 = P_\times,P_4 = P_<$ 
are predicates (of $\tau_1$) with 
1,1,3,3,2 places respectively such that using the definitions in 
clauses (B)(a),(b),(c) below, $P^M_\ell$ are 
defined in clauses (B)(d) below
\end{enumerate}
\sn
\item[(B)]
\begin{enumerate}
\item[(a)]  for $a \in M,\lev(a,M)$ is the maximal $k$ such that there are
$a_0 <^M_1 a^M_1 < \ldots <^M_1 a_k = a$; so necessarily $a_0 = a_*$ 
\sn
\item[(b)]  $\height(M) = \max\{\lev(a,M):a \in M\}$
\sn
\item[(c)]   for $k < \height(M)$ let 
$\iota = \iota(k,M) \in \{1,2\}$ be such that if $b$
 is a $<^M_2$-immediate successor of $a$ and $k = \lev(a,M)$ \then \,
$G \models \varphi_\iota[a,b]$, in the unlikely case both
$\iota=1$ and $\iota=2$ 
are as required we use $\iota=1$
\sn
\item[(d)]
\begin{enumerate}
\item[$(\alpha)$]  $P^M_0 = \{a_*\}$
\sn
\item[$(\beta)$]  $P^M_1 = \{a \in M:\lev(a,M)=1\}$
\sn
\item[$(\gamma)$]  $P^M_2 = \{(a,b,c):a,b,c \in M$
  and $\bbN \models ``\lev(a,M) + \lev(b,M) = \lev(c,M)"\}$
\sn
\item[$(\delta)$]  $P^M_3 = \{(a,b,c):a,b,c \in M$
  and $\bbN \models ``\lev(a,M) \times \lev(b,M) = \lev(c,M)"\}$
\sn
\item[$(\varp)$]  $P^M_4 = \{(a,b):N 
\models ``\lev(a,M) < \lev(a,N)"\}$.
\end{enumerate}
\end{enumerate}
\end{enumerate}
\end{definition}

\begin{definition}
\label{b9}
1) Let $\iota \in \{1,2\}$.

We say $N$ is the $\iota$-successor of $M$ in $\bold M_{G,a}$ \when
\, for some $k$
\mn
\begin{enumerate}
\item[(a)]  $M,N \in \bold M_{G,a}$ so $G \in K_n$ for some $n$
\sn
\item[(b)]  $M \subseteq N$ as models so $M = N \rest |M|$,
  recalling $|M|$ is the set of elements of $M$
\sn
\item[(c)]  $k = \height(M)$ and height$(N) = k +1$
\sn
\item[(d)]  $b \in N \backslash M$ \Iff \, $\lev(b,N) = k+1$ \Iff \,
  $b \in G \backslash M$ and for some\footnote{No real harm to demand here
(and in \ref{b4}) "unique"} 
$a \in M$ we have $\lev(a,M) = k$ and $G \models \varphi_\iota[a,b]$.
\end{enumerate}
\mn
2) We may omit $\iota$ above \when \,: $\iota=1$ \Iff \, $(*)$ holds where:
\mn
\begin{enumerate}
\item[$(*)$]  there is $c \in M$ such that $\lev(c,M) = k = 
\height(G)$ and for some $b \in G \backslash M$ we have $G
 \models \varphi_0[b,c]$; yes not $\varphi_1!$ but for case B there is
 no difference. 
\end{enumerate}
\mn
3) For a sentence $\psi$ in the vocabulary $\tau_1 \cup \tau_*$, (see
\ref{b2}, \ref{b4}), for $M,N \in \bold M_{G,a}$, we 
say $N$ is the $\psi$-successor of
$M$ \when \, for some $\iota \in \{1,2\},N$ is the 
$\iota$-successor of $M$ and $(G,M) \models \psi \Leftrightarrow (\iota=1)$.

\noindent
4) For $G \in K_n,a_* \in G$ and $M \in \bold M_{G,a_*}$ we define
$\bbN_{M,a_*}$ as the following structure $N$ with the vocabulary of
   number theory:
\mn
\begin{enumerate}
\item[$(a)$]  set of elements $\{a/E^M:a \in M\}$
\sn
\item[$(b)$]  $0^N = a_*/E^M = P^N_0$
\sn
\item[$(c)$]  $1^N = P^M_1$
\sn
\item[$(d)$]  if $\bold a_\ell = a_\ell/E^M \in N,a_\ell \in M$ for
$\ell = 1,2,3$ then
\sn
\begin{enumerate}
\item[$(\alpha)$]  $N \models ``\bold a_1 + \bold a_2 = \bold a_3"$
iff $(a_1,a_2,a_3) \in P^M_2$
\sn
\item[$(\beta)$]  $N \models ``\bold a_1 \times \bold a_2 = \bold a_3"$
iff $(a_1,a_2,a_3) \in P^M_3$
\sn
\item[$(\gamma)$]  $N \models ``\bold a_1 < \bold a_2"$
iff $(a_1,a_2) \in P^M_4$.
\end{enumerate}
\end{enumerate}
\end{definition}

\begin{claim}
\label{b12}
1) If $\iota,G,M,a$ are as in \ref{b9}(1) \then \, there is at most one
$\iota$-successor $N$ of $M$ in $\bold M_{G,a}$.

\noindent
1A) For some $\psi_* \in \bbL(\tau_1 \cup \tau_*)$, 
being a $\psi_*$-successor is equivalent to
being a successor.

\noindent
2) For a given $G \in K_n,a \in G,M \in \bold M_{G,a}$ and $\psi \in 
\bbL(\tau_1 \cup \tau_*)$ 
there 
is at most one $\psi$-successor $N$ 
of $M$ in $\bold M_{G,a}$.

\noindent
3) For $G \in K_n$ and $a \in G$
   there is one and only one sequence $\langle M_k = 
M_{k,a}:k \le k_{G,a}\rangle$ such that:
\mn
\begin{enumerate}
\item[$(a)$]  $M_k \in \bold M_{G,a}$
\sn
\item[$(b)$]  $M_0$ has universe $\{a\}$
\sn
\item[$(c)$]  $M_{k+1}$ is the successor of $M_k$ in $\bold M_{G,a}$,
  recall \ref{b9}(2)
\sn
\item[$(d)$]  if $k=k_{G,a}$ \then \, there is no $N \in 
\bold M_{G,a}$ which is the successor of $M_k$ in $\bold M_{G,a}$
\end{enumerate}
\mn
3A) Above $\bbN_{M_k,a}$ is isomorphic to $\bbN \rest \{0,\dotsc,k\}$.  Also
for every sentence $\psi \in \bbL$ or even $\in \bbL_{\LFP}$ 
in the vocabulary of number theory there is
a sentence $\varphi \in \bbL_{\LFP}(\tau_*)$ such that $\bbN_{M_k,a}
\models \psi \Rightarrow M_{k,a} \models \varphi$; of course, $\varphi$
depends on $\psi$ but not on $G,a$ (and $k$).

\noindent
4) In the LFP logic for $\tau_*$, we can find a sequence
$\bar\varphi$ of formulas with\footnote{The variable $y$ stands for
  the parameters $a$; instead we may define in \ref{b9} one model
  $M_k$ coding all $M_{a,\ell} \in \bold M_{G,a}$ for $\ell \le k,a
  \in G$.} variable $x_0,\dotsc,y$ such that: for any 
$G \in K_n$ and $a \in G$, the sequence $\bar\varphi$ substituting $y$
by $a$ defines $M = M_{G,a}$ which is $M_{k,a}$ 
for $k=k_{G,a,\psi}$ from part (3).

\noindent
4A)  For $\psi$ as in \ref{b9}(3), i.e. $\psi \in \bbL(\tau_1 \cup
\tau_*)$ , a sentence, the parallel of \ref{b12}(3),(4),(5) 
 holds for ``$\psi$-successor", (so we should write 
$M_{G,a,\psi}$ instead $ M_{G,a}$).

\noindent
5) Letting $\height(G) = \max\{\height(M_{G,a}):a \in G\}$, in LFP 
logic there is $\varphi_*(x)$ such that $G \models
\varphi_*(a)$ \Iff \, $a \in G$ and $\height(M_{G,a}) = \height(G)$ \Iff \, 
 for every $a_1 \in G,\height(M_{G,a_1}) \le \height(M_{G,a})$. 

\noindent
6) For any sentence $\psi$ in the vocabulary of number theory (in
first order or LFP logic) there is a sentence $\varphi$ in
induction logic for $\tau_*$ (recallng \ref{b2})
such that for any $G,G \models \varphi$
iff $\bbN \rest \{0,\dotsc,\height(G)\} \models \psi$.
\end{claim}

\begin{PROOF}{\ref{b12}}
1) Read the Definition \ref{b9}(1).

\noindent
1A) Read \ref{b9}(2).

\noindent
2) Read \ref{b9}(3) and recall part (1).

\noindent
3) We choose $M_k$ and prove its uniqueness by induction on $k$ till we
are stuck.  Recalling \ref{b12}(3)(d) we are done.

\noindent
3A) Easy.

\noindent
4),4A) Should be clear.

\noindent
5) We can express by induction when ``$\lev(b_1,M_{G,a_1}) \le
\lev(b_2,M_{G,a_2})"$.

\noindent
6) Should be clear but we elaborate.

Recall the formula $\varphi_*(x) \in \bbL_{\LFP}(\tau_*)$ from
\ref{b12}(5).  By the choice of $\varphi_*$ necessarily 
for some $a_*,G_n \models \varphi_*[a_*]$
(as in a finite non-empty set of natural numbers there is a maximal member) so
$\height(a_*,G_n) = \height(G_n)$.

Now for a given $\psi$, let $\varphi \in \bbL_{\LFP}(\tau_*)$
say: for some (equivalently every) $a \in G_n$ such that $G_n \models
\varphi_*(a)$, the model $\bbN_{G_n,a_n}$ defined in \ref{b9}(4),
which is isomorphic to $\bbN \rest \{0,\dotsc,\height(a,M_{G_n,a_n})\}$, see
\ref{b12}(3A), satisfies $\psi$.
\end{PROOF}

\begin{conclusion}
\label{b13}
We have ``$G_n$ fail the 0-1 law for the LFP logic, moreover 
for some $\varphi \in \bbL_{\ind}(\tau_*)$ we have
$\Prob(G_n \models \varphi)$ has $\lim - \sup= 1$ and 
$\lim \inf =0$".
\end{conclusion}

\begin{PROOF}{\ref{b13}}
Should be clear by the above, in particular \ref{b12}(6), 
see \ref{d3}, \ref{d4}(2) for details on the probabilistic estimate
needed for \ref{a7}(2) on which we rely.
But we elaborate.

Note that just the following is not sufficient:
\mn
\begin{enumerate}
\item[$(*)_1$]  some $\bar\varphi,\bold m$ satisfies
\sn
\begin{enumerate}
\item[$(a)$]  $\bar\varphi$ an interpretation scheme, see \ref{z3d}
\sn
\item[$(b)$]  $\bold m$ is a function from the class of finite graphs
to $\bbN$, depending on the isomorphism type only
\sn
\item[$(c)$]  for every $m$ for every random enough $G_n,\bold m(G_n)
\ge m$
\sn
\item[$(d)$]  for random enough $G_n,\bar\varphi$ defines an
isomorphic copy of $\bbN \rest \{0,\dotsc,\bold m(G_n)\}$.
\end{enumerate}
\end{enumerate}
\mn
However it is enough if we add, e.g.
\mn
\begin{enumerate}
\item[$(*)_2$]  $\bold m(G_n) \ge \log_2(\log_2(n))$.
\end{enumerate}
\mn
\underline{Why it is enough}?  Let $\psi$ be a first order sentence in
the vocabulary such that $\bbN \rest \{0,\dotsc,k\} \models \varphi$
\Iff \, $\log_*(k)$ belong to $\{10n + \ell:\ell = \{0,1,2,3,4\}$ and $n \in
\bbN\}$.

Now use the interpretation $\bar\varphi$, i.e. we use \ref{b12}(6) in
our case.
\medskip

\noindent
\underline{Why $(*)_1 + (*)_2$ holds}:  By \ref{d3}, \ref{d4}(1) and
\ref{b12}(3A). 
\end{PROOF}
\newpage

\section {Revisiting induction}

As discussed in \S2, we need to prove that for random enough
$G_n,\height(G_n)$ is large enough, equivalently, for some $a \in G_n$
(we shall prove that even for most), 
$\height(a,G_n)$ is large enough.  For this a more
detailed specific statement is proved - see $(*)$ in the proof of
\ref{d3}.  That is, we prove that for most $a \in G_n$ (for
random enough $G_n$): on the one hand $M_{G,a,=k}$ is not too large,
and, on the other hand, is not empty; and for Case A, even not too small.  
The computation naturally depends on what
$\eta_{G,a}$ is, see \ref{b17}(3), this is a delicate point.  

For Case B, things are simpler.  For each $k$ we ask if there is an
$R_1$-edge out of $M_{G,a,=k}$ to $G \backslash M_{G,a,k}$.  If
there is, clearly $M_{G,a,=k+1}$ will be quite small but not empty.
If not, then necessarily $M_{G,a,=k}$ has $\le n^{\alpha^*_2-\zeta}$ members 
hence the number of
$R_2$-neighbors of members of $M_{G,a,=k}$ cannot be too
large (well $< n^{\alpha^*_2} n^{\alpha^*_2}$) so we are done.

Case A seems harder, so we simplify considering only small enough $k$,
see \ref{d4}, hence we can consider all possible $\eta$'s of length
$k$, that is, summing the probability of the ``undesirable" events on
all of them; so if each has small enough probability, even the unions 
of all those events has small enough probability.  Now we divide the
$\eta$'s to those which are ``reasonable candidates to be
$\eta_{G_n,a}$" and those which are not.  For the former $\eta$'s, for
almost all $a \in G_n$ there is a pre-$(\eta,0,k_*)$-path starting
with $a$.  So it is enough to prove that for almost all $a \in
G_n,\eta_{G,a} \rest k^*$ is one of those former $\eta$'s where $k_*$
is the relevant large enough height, e.g. $\ge \lfloor
\log_2(\log_2(n)) \rfloor$.  For this
it is enough to prove that the other $\eta$'s cannot occur and this
is what we do.

We in this section fulfill promises from \S2 (and \S1).
\begin{context}
\label{d1}
As in \ref{b2}.
\end{context}

\noindent
Below we shall use
\begin{definition}
\label{b17}
1) For $G \in K_n$ we define $M_k(a,G) = M_{G,a,k}$ by induction as in
\ref{b12}(3) for $\psi = \psi_*$ from \ref{b12}(1A) and also
$k=k_{G,a}$ as there and $\height(G)$ as in \ref{b12}(6).

\noindent
2) Let $M_{G,a,=k} = M_{G,a,k} \backslash \cup\{M_{G,a,m}:m < k\}$ and
similarly $M_{G,a,<k}$.

\noindent
3) Let $\eta = \eta_{G,a}$ be the following sequence of length
   $k_{G,a}$: if $\ell < k_{G,\ell}$, then $\eta(\ell) =
   \iota(\ell,M_{G,a}) \in \{1,2\}$ from Definition \ref{b4}(B)(c).
\end{definition}

\begin{claim}
\label{d3}
For small enough $\varepsilon \in (0,1)_{\bbR}$, 
for random enough $G_n$, for some $a \in [n],k_{G_n,a} 
\ge k^* = \lfloor n^\varepsilon \rfloor$ in case B and $k_{G_n,a} \ge
\lfloor \log(\log(n)) \rfloor$ in case A.
\end{claim}

\begin{remark}
\label{d4}
It would be nice to use an $\eta \in {}^{\bbN}\{1,2\}$ defined
similarly to \ref{a4}, say such that $\gamma_n \in
[\alpha^*_0,\alpha^*_2 + \alpha^*_2]$, but this is not clear.
In case B, in the proof the problem is that the $\gamma_n$-s 
from \ref{a4} may be  very near to $\alpha^*_0$.
Also the parallel problem for case A is that 
the answer to the question asked there is near the 
critical stage, so we are
not almost sure about the answer.
\end{remark}

\begin{PROOF}{\ref{d3}}
For case A,  we presently prove it, e.g. for 
$k^* = \lfloor \log_2(\log_2 n) \rfloor$ and for Case B, $\lfloor
n^\varp \rfloor$, an overkill,
but 
this suffices for the failure of the 0-1 law.  We intend to fill
the general case elsewhere.
Actually for any $\varepsilon \in (0,1)_{\bbR}$ we 
can get $k^* = \lfloor n^{1-\varepsilon} \rfloor$.

Let $\zeta \in (0,1)_{\bbR}$ be small enough and $k^*$ be as above.

Clearly it is enough to prove:
\mn
\begin{enumerate}
\item[$(*)$]  if $a \in [n]$ and $k < k^*$ \then \, the probability
  that at least one of the following $(i)_{a,k},(ii)_{a,k},(iii)_{a,k}$ fails
  (assuming $\ell < k \Rightarrow (i)_{a,\ell} \wedge
  (ii)_{a,\ell})$ is small enough; $< \frac{1}{k \log(n)}$ suffices,
  being $< \frac{1}{k n^i}$ for each $i$ for large enough $n$ is natural)
\sn
\item[${{}}$]  $(i)_{a,k} \quad k \le k_{G_n,a}$
\sn
\item[${{}}$]  $(ii)_{a,k} \quad M_{G_n,a=k}$ has 
$\le n^{\alpha^*_2 + \alpha^*_2}$ elements
\sn
\item[${{}}$]  $(iii)_{a,k} \quad M_{G_n,a,k}$, (noting that
  $(i)_{a,k}$ implies $M_{G_n,\ell,k} \ne \emptyset$) has 

\hskip25pt $\ge n^{\alpha^*_0 -
\zeta}$ elements \footnote{We can use $\ge n^{\alpha^*_2 - 
\alpha^*_0 -\zeta}$}  except when $k=0$, \underline{not need} for 
Case B.
\end{enumerate}
\mn
Why does $(*)$ hold?  
\medskip

\noindent
\underline{Case 1} Case B of the Context \ref{b2} and Definition \ref{z13}

We are given $n \ge 1$ and $a \in [n]$; we draw the edges in $k$ stages 
so by induction on $k$.  For $k=0$ draw the edges starting with $a$
(of both kinds, an overkill), i.e. for $\iota \in \{1,2\}$
the truth value of $R_\iota(a,b)$ for every $b \in
[n] \backslash \{a\}$, hence we can compute $M_1(a,G)$.

The induction hypothesis on stage $k$ is that $\langle M_{G,a,i}:i \le
k\rangle$ have been computed and we have drawn the truth value of
$R_\iota(c,b)$ for $b \in \cup\{M_{G,a,i}:i < k\}$ and $c \in [n]
\backslash \{b\}$.  If $k < k_*$ we now draw the edges $R_\iota(b,c)$ for
$b \in M_{G,a,=k}$ and any $c \ne b$; 
actually the $c \in M_{G,i,k}$ are irrelevant and so we can compute
$M_{G,a,k+1}$.  Now we ask: if $(i)_{a,m} + (ii)_{a,m}$ holds
for $m \le k$ what is the probability that $(i)_{a,k+1} +
(ii)_{a,k+1}$? (recalling $(iii)_{a,k+1}$ is irrelevant), 
i.e., is it small enough? 
This is easy and as required.
\bigskip

\noindent
In details, we ask
\smallskip

\noindent
\underline{Question}:  Are there $c \in [n] \backslash M_{G,a,k}$ and
$b \in M_{G,a,=k}$ such that $(b,c) \in R^G_1$?

First note
\mn
\begin{enumerate}
\item[$(*)$]  if 
 $|M_{G,a,=k}| \ge n^{\alpha^*_2-\zeta}$ 
then the probability that the answer is no is $< 1/2^n$.
\end{enumerate}
\mn
[Why?  We have $M_{G,a,=k} \times ([n] \backslash M_{G,a,k})$
independent drawings so their number is $\ge n^{\alpha^*_2-\zeta} n/2$, 
each with probability $\frac{1}{n^{1+\alpha^*_1}}$
of success and $(1 + \alpha^*_2-\zeta) - (1 + \alpha^*_1) =
\alpha^*_2-\zeta - \alpha^*_1 >0$ so the probability of the no answer
is $(1 - \frac{1}{n^{1+\alpha^*_1}})^{n^{(1 +\alpha^*_2-\zeta)}} \sim
1/e^{(n^{\alpha^*_2-\zeta - \alpha^*_1)/2})}$; clearly more than enough.] 

By $(*)$ it suffices to deal with the following two possibilities.
\medskip

\noindent
\underline{Possibility 1}:  The answer is yes.

In this case $M_{G,a,k+1}$ is well defined and
$\iota(k,M_{G,a,k})=1,M_{G,a,=k+1} = \{c:c \in G \backslash 
M_{G,a,k}$ and $(b,c) \in R^G_1$ for some $b \in M_{G,a=k}\}$.  Now for
each $c \in [n] \backslash M_{G,a,k}$ and $b \in M_{G,a,=k}$ the
probability of $(c,b) \in R^G_1$ is $\frac{1}{n^{1+\alpha^*_1}}$ 
hence by the independence of the drawing, recalling
$|M_{G,a,=k}| \le n^{\alpha^*_2 + \alpha^*_2}$ the probability of
$|M_{G,a,=k+1}| \ge n^{\alpha^*_2 +\alpha^*_2}$ is negligible, e.g. 
$< 2^n$ so can be ignored.
Also by the possibility we are in, $M_{G,a,= k+1} \ne \emptyset$. 
\medskip

\noindent
\underline{Possibility 2}:  The answer is no and 
$|M_{G,a,=k}| \le n^{\alpha^*_2-\zeta}$.

This is easy, too, recalling that 
almost surely for every $a' \in G$ the number of $R_2$-neighbors 
in the interval $[n^{\alpha^*_2} -
n^{\alpha^*_2(1-\zeta)},n^{\alpha^*_2} + n^{\alpha^*_2(1-\zeta)}]$
and so the probability that $M_{G,a,=k+1}$ is too large is negligible.
\medskip

\noindent
\underline{Case 2}:  Case A of the context \ref{b2}

Here it helps to use ``$\varphi_0,\varphi_1$ are distinct".

Now it suffices to prove:
\mn
\begin{enumerate}
\item[$(*)_1$]  for random enough $G_n$, for\footnote{We allow few
$a$'s for which this fails.  It suffices to have ``for some $a$",
this helps for larger $k$.} $a \in M$, the following has
negligible probability of failure: $(i)_{a,k},(ii)_{a,k},(iii)_{a,k}$.
\end{enumerate}
\mn
Note that for this it seems more transparent  
to\footnote{As then we can consider all the relevant sequences
$ \eta $, (and more)} 
assume $k <
\log_2(\log_2(n))$ and to translate $(*)_1$ to statement on paths.
\mn
\begin{enumerate}
\item[$(*)_2$]  For $k \le k_*$ let
\sn
\begin{enumerate}
\item[$(a)$]  $\Omega_{=k}$ be the set of $\eta \in {}^k\{1,2\}$ such
that $\alpha(\eta) := |\eta^{-1}\{2\}|
\cdot \alpha^*_2 - |\eta^{-1}\{1\}| \cdot \alpha^*_1$ belongs to the
interval $[0,\alpha^*_2 + \alpha^*_2 + \zeta]$
\sn
\item[$(b)$]  $\Omega_k$ be the set of $\eta \in {}^k\{1,2\}$ such
that for every $\ell < n$ the sequence $\eta \rest \ell = \langle
\eta(0),\dotsc,\eta(\ell-1)\rangle$ belongs to $\Omega_{=\ell}$
\sn
\item[$(c)$]  $\Delta_k$ be the set of $\eta \in {}^k\{1,2\}$ such
that $\eta \in \Omega_{=k}$ and even $\eta \in \Omega_k$ 
but $\alpha(\eta) \le \alpha^*_0 - \zeta$
\end{enumerate}
\sn
\item[$(*)_3$]  recalling that $\zeta \in (0,1)_{\bbN}$ is small
enough, for any random enough $G_n$, for every $a \in M$ the following
 has probability $\le 1/n^\zeta$: 
\sn 
\begin{enumerate}
\item[$\bullet$]  for some $k \le \lfloor \log(\log(n)) \rfloor$ and
  $\eta \in {}^k\{1,2\}$ at least one of the following holds
\sn 
\item[${{}}$]  $(a)_\eta \quad \eta \in \Omega_k$ but there is no 
pre-$(\eta,0,k)$-path in $G_n$ starting with $a$
\sn
\item[${{}}$]  $(b)_\eta \quad \eta \in \Delta_k$ but there is a
pre-$(\eta,0,k+1)$-path in $G_n$ from $a$ 

\hskip25pt to some $b \in G_n
\backslash \{a\}$ such that $G_n \models (\exists x)\varphi_0(b,x)$.
\end{enumerate}
\end{enumerate}
\mn
Otherwise the proof is as in the earlier case.
\end{PROOF}

\begin{claim}
\label{d14}
Let $G_n = G_{n;\bar\alpha^*}$, i.e. we are in Case B.  For some
$(\tau_{\bbN},\tau_{\dn})$-scheme $\bar\psi$, for every random
enough $G_{n;\bar\alpha^*,\psi}$ defines in $G_n$ a structure isomorphic to
$\bbN_{<n}$. 
\end{claim}

\begin{PROOF}{\ref{d14}}
We make a minor change in Definition \ref{b9}(1),(d).
\medskip

\noindent
\underline{Clause $(d)^+$}:  We require the $a \in M$ is unique.  

This makes no real difference above because the probability of the
occurance if even one ``$b$
with two predecessors" is small and we just need that there is one;
this makes $M_{G,a,k}$ smaller but not empty.

We start as in the proof of \ref{d3}, for $k_* = \lfloor n^\zeta
\rfloor$, we use Case 1 but for stage $k$ we draw the truth values of
$\bbR_\iota(b,c)$ only when $b \in M_{G,a,=k}$ and $c \in [n]$ but $c
\notin M_{G,a,0} \cup \ldots M_{G,a,k-1}$.

So there is $b \in M_{G,a,=k_*}$ and there is a unique sequence
$\langle a_\ell:\ell \le k_*\rangle,a_0 = a,a_{k_*}=b$ and
$(a_\ell,a_{\ell +1}) \in R^{G_n}_\iota$ where $\iota$ is such that
$M_{G,a,\ell+1}$ is the $\iota$-successor of $M_{G,a,\ell}$
and so there are formulas $\psi_2 \in \bbL_{\LFP}(\tau_{\dg})$ 
such that not depending on
the pair $(a,n),\psi_2(G,a,b) = \{(a_\ell,a_i):\ell \le i \le k_*\}$.

Now
\mn
\begin{enumerate}
\item[$\bullet$]  the probability of the following event is negligible
  ($c<2$): for some $d_1 \ne d_2 \in [n] \backslash M_{n,k_*}$ for
  every $c$:  if $\psi_2(c,c,a,b)$ then $(d_1,c) \in R^G_2
  \Leftrightarrow (d_2,c) \in R^G_2$.
\end{enumerate}
\mn
Ignoring this event, the following formulas defines a linear order on
$[n] \backslash M_{G,a,k_*}$:
\mn
\begin{enumerate}
\item[$\bullet$]  $\psi_3(d_1,d_2,a,b)$ say: for some $c$ we have
$\psi_2(c,c,a,b) \wedge R_2(d_2,c) \wedge \neg R_2(d_1,c)$ and for
any $c_1$ if $\psi_2(c_1,c,a,b)$ then $R_2(d_1,c') \leftrightarrow
R_2(d_2,c')$.
\end{enumerate}
\mn
So $(x,y,a,b)$ defines a linear order on $[n] \backslash
M_{G,a,k_*}$ which has $\ge n - \lfloor n^\zeta \rfloor$ 
elements.  We get using the same trick $\psi_4 \in 
\bbL_{\LFP}(\tau_{\dn})$ and $\psi_4(x,y,a,b)$
defines a linear order on $[n]$.  Now the formulas
$\psi_0,\dotsc,\psi_4$ do not depend on $n$.  Also for some
$\psi_5 \in \bbL_{\LFP}(\tau_{\dg})$
\mn
\begin{enumerate}
\item[$\bullet$]  $G \models \psi_5(a,b)$ iff
$\psi_4(-,-,a,b)$ defines a linear order (on $G$)
\newline
and for some $\psi_0 \in \bbL_{\LFP}(\tau_{\dg})$
\sn
\begin{enumerate}
\item[$\odot$]  $G \models \psi_0[c_1,a_1,b_1,c_2,a_2,b_2]$ iff:
\sn
\item[${{}}$]  $(a) \quad (a_1,b_1),(a_2,b_2) \in \varphi_5(G)$
\sn
\item[${{}}$]  $(b) \quad |\{c:G \models \psi_4[c,c_1,a_2,b_1)\}| =
  |\{c:G \models \psi_4[c,c_2,a_2,b_2]\}|$.
\end{enumerate}
\end{enumerate}
\mn
So the interpretation should be clear.
\end{PROOF}

\begin{conclusion}
\label{d17}
[Case B]  For some $\psi \in \bbL_{\LFP}(\tau_{\dg})$ for every random
enough $G_n = G_{n;\bar\alpha^*}$ we have:
\mn
\begin{enumerate}
\item[$\bullet$]  $G_n \models \psi$ iff $n$ is even.
\end{enumerate}
\end{conclusion}
\newpage

\bibliographystyle{alphacolon}
\bibliography{lista,listb,listx,listf,liste,listz}

\end{document}